\documentclass[11pt]{amsart}
\usepackage{amsfonts}
\usepackage{amssymb}
\usepackage[dvips]{graphics,color}

\voffset=-.2in \textheight=8.75in \theoremstyle{plain}
\newtheorem{cor}{Corollary}
\newtheorem{thm}{Theorem}

\begin{document}

\title[Combinatorial Identities]
{\bf On Proofs of Certain Combinatorial Identities}

\author[Grossman, Tefera and Zeleke]
{George Grossman, Akalu Tefera and  Aklilu Zeleke }
\address{ Department of Mathematics \\
Central Michigan University \\
Mt. Pleasant, MI 48859 \\
gross1gw@cmich.edu}
\address{ Department of Mathematics \\
Grand Valley State University \\
Allendale , MI 49401 \\
TeferaA@gvsu.edu}
\address{ Department of Mathematics \\
Alma College \\
Alma, MI 48801 \\
zeleke@alma.edu}


\maketitle

\section*{Abstract}

In this paper we formulate combinatorial identities that give
representation of positive integers as linear combination of even
powers of $2$ with binomial coefficients.  We present side by side
combinatorial as well as computer generated proofs using the
Wilf-Zeilberger(WZ) method.

\section{Introduction}

It is known that every integer can be written as a sum of integral
powers of 2. A somewhat related problem is to find for every
positive integer $n$ a positive integer $k$ depending on $n$ with
$k(n) < k( n+1)$ and  integer coefficients $a_{i}, i = 0,\, 1,\,
\ldots,\, k-1$ such that
\begin{equation}
\label{eq:1}
 n = \sum_{i=0}^{k-1}a_{i}\,2^{2\,i}\ .
\end{equation}
The background and motivation for this problem lies in studying
the zeros of the $j-th$ order polynomial of the generalized
Fibonacci sequence given by
\begin{equation}
\label{eq:2}
 F_{j}(x)= x^{j} - x^{j-1} - \cdots - x - 1 \ .
\end{equation}
 For studies related to the
positive zeros of (2) we refer the reader to the papers by Dubeau
(\cite{D89}) and Flores (\cite{F67}). It can be shown (see
\cite{GN99})  that for $j$ even
\begin{equation}
F_{j}(x)=
(x-2+\epsilon_{j})(x+1-\delta_{j})(x^{j-2}+a_{j-3}x^{j-3}+...+a_{1}x+a_{0})
\ ,
\end{equation}
where $-1+\delta_{j}$ and $2-\epsilon_{j}$ are  the negative and
positive zeros of (\ref{eq:2}). Here  $\{\delta_{j}\}$ and
$\{\epsilon_{j}\}$ are positive, decreasing sequences. In a recent
paper, Grossman and Zeleke (\cite{GZ03}) have found an explicit
form for the $a_{i}$'s in terms of $\epsilon_{j}$ and $\delta_{j}$
for $j \geq 4$. The explicit expressions for $a_{i}$ as well as
special cases led to some interesting identities. In this paper we
present different proofs of three such identities that are
hypergeometric.  The paper is organized as follows. In section 2,
we formulate the main results. In section 3 we provide
combinatorial proofs. This requires first finding combinatorial
interpretations by counting words of certain properties and
defining an appropriate sign reversing involution which we call
``involution''. Gessel and Stanley discuss the mathematical theory
related to such proofs in (\cite{GS95}).
 In section 4, we present computer
generated proofs of the main results. It is to be noted that there
are philosophical arguments over computer-based proofs to
mathematical proofs in general. It will be clear from sections 3
and 4 that the WZ method gives a unified and structured approach
to proving identities of this type. Introductions to the WZ method
can be found among others in the book \textit{A = B}
(\cite{PWZ96}) or the website
\texttt{http://mathworld.wolfram.com/Wilf-ZeilbergerPair.html}.
Throughout this paper we denote the set $\{k, k + 1, k + 2, \dots
\}$ for $k \in \mathbb{Z}$ by $\mathbb{N}_{k}$.

\section{Main Results}

\begin{thm} \cite{GZ03}
\label{Thm1}
\[
 \sum_{k = 0}^{n} (-1)^{n + k}{{n + k + 1} \choose {2\,k+
1}}2^{2\,k}=  n + 1, \;\; n \in \mathbb{N}_{0} \ .
\]
\end{thm}

\textbf{Remark:} The following theorems  will show that the
coefficients $a_{i}$'s in the  expansion of positive integers are
not unique.

\begin{thm}
\label{Thm2}
\[
 \sum_{k = 0}^{n + 1}(-1)^{n + k + 1} {{n + k + 1}
\choose {2\,k}}2^{2\,k}=  2\,n + 3, \;\; n \in \mathbb{N}_{-1}
 \ .
\]

\end{thm}

\begin{thm} \cite{GZ03}
\label{Thm3}
\[
\sum_{k = 0}^{n - 1}\sum_{m = 2\,k + 2}^{n + 1 + k} (-1)^{m + k}
{{n + k + 1} \choose {m}}2^{m - 1} = n\,(n + 1),\;\; n \in
\mathbb{N}_{1} \ .
\]
\end{thm}

\bigskip

\section{Proofs of the Main Results}

\subsection{Combinatorial Proofs}

For a fixed  $n \in \mathbb{N}_{k}$, $k \in \mathbb{Z}$, consider
the set $S$ of words in the alphabet $\{a,\, b,\, c\}$ such that
\[
2\,(\mbox{the number of a's}) + 1\,(\mbox{the number of b's}) +
1\,(\mbox{the number of c's}) = p(n),
\]
where $p(n)$ is some polynomial of $n$. For $w \in S$, define the
weight $Wt(w)$ of $w$ by
\[Wt(w) = (-1)^{(\mbox{the number of
a's})}.\]

\textbf{Proof of Theorem \ref{Thm1}.}

For $n \in \mathbb{N}_{0}$, take $p(n) = 2\,n + 1$. Then
\begin{eqnarray*}
\sum_{w \in S}Wt(w) & = & \sum_{k = 0}^{n} (\mbox{the number of
w's with } (n - k)\,\,  a's)\,Wt(w)\\
& = & \sum_{k = 0}^{n} {{n + k + 1} \choose {2\,k + 1}} 2^{2\,k +
1}\,(-1)^{n + k}\\
& = & 2\,\sum_{k = 0}^{n} (-1)^{n + k}\,{{n + k + 1} \choose {2\,k
+ 1}} 2^{2\,k} \ .
\end{eqnarray*}

Consider now $S  = T \cup (S - T)$, where $T$ is the set of all
words in $S$ of the form $c^{m}\,b^{2\,n + 1 - m}$ for $m =
0,\,1,\,\ldots, 2\,n + 1$. Here by, for $x \in \{a,\,b,\,c\}$ the
notation $x^{0}$ is used to denote the empty word. Define an
``involution'' as follows:

For $w \in S$ read $w$ left to right until you either get  an $a$,
or $b\,c$. If it is an $a$, make it a $b\,c$. If it is a $b\,c$,
make it an $a$. This changes the sign of $Wt(w)$ and is an
involution. Note that $T$ has $2\,n +  2$ elements each of weight
$(-1)^{0} = 1$. From the involution it is clear that the sum of
the weights of the elements of $S - T$ is $0$. Thus $\displaystyle
\sum_{w \in S}Wt(w) = 2\,n + 2$. Hence the theorem follows.

\newpage

\textbf{Proof of Theorem \ref{Thm2}.}

Let $n \in \mathbb{N}_{-1}$ and $p(n) = 2\,n + 2$. Then
\begin{eqnarray*}
\sum_{w \in S}Wt(w) & = & \sum_{k = 0}^{n + 1} (\mbox{the number
of
w's with } (n + 1 - k)\,\,  a's)\,Wt(w)\\
& = & \sum_{k = 0}^{n} {{n + k + 1} \choose {2\,k}} 2^{2\,k}\,(-1)^{n + k + 1}\\
& = & \sum_{k = 0}^{n} (-1)^{n + k + 1}\,{{n + k + 1} \choose
{2\,k}} 2^{2\,k} \ .
\end{eqnarray*}

Partition $S$ as in the proof of Theorem 1 with $T$ the set of all
words in $S$ of the form $c^{m}\,b^{2\,n + 2 - m}$ for $m =
0,\,1,\,\ldots, 2\,n + 2$. $T$ has $2\,n + 3$ elements each of
weight $(-1)^{0} = 1$.

Using the same ``involution'' as in the proof of Theorem 1 the sum
of the weights of the elements of $S - T$ would be $0$ and hence
$\displaystyle \sum_{w \in S}Wt(w) = 2\,n + 3$.

\textbf{Proof of Theorem \ref{Thm3}.}

For $n \in \mathbb{N}_{1}$, consider the set $S$ of words in  the
alphabet $\{a,\, b,\, c\}$ such that
\[
1\,(\mbox{the number of a's}) + 1\,(\mbox{the number of b's}) +
1\,(\mbox{the number of c's}) = n + 1 + k,
\]
for some $k \in \{0,\, \ldots, n - 1\}$ and $(\mbox{the number of
b's}) + (\mbox{the number of c's})$ is at least $2\,k + 2$. For $w
\in S$, define the weight $Wt(w)$ of $w$ by
\[Wt(w) = (-1)^{(\mbox{the number of
a's})}.\]

Then

\begin{eqnarray*}
\sum_{w \in S}Wt(w) & = & \sum_{k = 0}^{n - 1}\sum_{m = 2\,k +
2}^{n + k + 1} (\mbox{the number of
w's with } (n + 1 + k - m)\,\,  a's)\,Wt(w)\\
& = & \sum_{k = 0}^{n - 1}\sum_{m = 2\,k + 2}^{n + k + 1} {{n + k
+ 1} \choose {m}} 2^{m}\,(-1)^{n + 1 + m + k}\\
& = & 2\,(-1)^{n + 1} \sum_{k = 0}^{n - 1}\sum_{m = 2\,k + 2}^{n +
k + 1} {{n + k + 1} \choose {m}} 2^{m - 1}\,(-1)^{m + k}\ . \\
\end{eqnarray*}

Read a word $w \in S$ from left to right. Count the number $b$ and
$c$ until the sum is $3$. Thus $w$ has the form
$a^{l}\,x\,a^{n}\,y\,a^{m}\,z\,\ast$ where $l, m, n \in
\mathbb{N}_{0}$ and $x,\, y,\, z \in \{b,\,c\}$. For such words,
define a mapping
 $\sigma$ as follows:
\[
\sigma(w) = \left\{ \begin{array}{r@{\qquad:\quad}l}
a^{l}\,x\,a^{n + 1}\,y\,a^{m}\,z\,\ast &  \mbox{ if } n, m
\mbox{ have same parity and } n \neq 1 \ , \\

a^{l}\,x\,a^{n - 1}\,y\,a^{m}\,z\,\ast &  \mbox{ if } n, m
\mbox{ have same parity and } n = 1 \ , \\

a^{l}\,x\,a^{n - 1}\,y\,a^{m}\,z\,\ast &  \mbox{ if } n, m
\mbox{ have opposite parity and } n \neq 0 \ , \\
a^{l}\,x\,a^{n + 1}\,y\,a^{m}\,z\,\ast &   \mbox{ if } n, m \mbox{
have opposite parity and } n =0 \ .
\end{array}
\right .
\]

Clearly $\sigma$ is an ``involution''. This involution is not
defined for elements of $S$ of length $n + 1$ and the number of
$b$'s and $c$'s exactly 2. There are $\displaystyle 4\,{{n + 1}
\choose {2}}~=~ 2\,n(n + 1)$ such words each of weight $(-1)^{n +
1}$ and hence the theorem follows.

\subsection{The WZ Method Proofs}

\textbf{Proof of Theorem \ref{Thm1}.}

Let $\displaystyle F(n,\, k) = {{n + k + 1} \choose {2\,k + 1
}}\frac{2^{2\,k}\,(-1)^{k + n + 1}}{n + 1}$ and let $\displaystyle
S(n) = \sum_{k = 0}^{n} F(n,\,k)$. We want to show that $S(n) = 1$
for all $n \in \mathbb{N}_{0}$. $F$ satisfies the recurrence
equation:\footnote{The recurrence equation is automatically
generated by a MAPLE package \texttt{EKHAD} which is available
from \texttt{htpp//www.math.rutgers.edu/\~{}zeilberg/}}
\begin{eqnarray}
\label{Thm1eqn1}
 F(n + 1,\, k) + F(n,\, k) =
G(n,\, k + 1) - G(n,\, k),
\end{eqnarray}
where $G(n,\, k) = R(n,\,k)\,F(n,\, k)$ and $\displaystyle
R(n,\,k) =  -\frac {k\,(2\,k + 1)}{(n + 1 - k)\,(n + 2)}$ \ .

By summing both sides of equation (\ref{Thm1eqn1}) with respect to
$k$ we get $S(n + 1)  - S(n) = 0$. Moreover,  $S(0) = 1$ and hence
$S(n) = 1$ for all $n \in \mathbb{N}_{0}$.

\bigskip

\textbf{Proof of Theorem \ref{Thm2}.}

Let $\displaystyle F(n,\,k) = {{n + k + 1} \choose
{2\,k}}\frac{2^{2\,k}\,(-1)^{k + n + 1}}{2\,n + 3}$ and let
$\displaystyle S(n) = \sum_{k = 0}^{n + 1} F(n,\,k)$. We want to
show that $S(n) = 1$ for all $n \in \mathbb{N}_{-1}$. $F$
satisfies the recurrence equation:$^1$
\begin{eqnarray}
\label{Thm2eqn1}
 F(n + 1,\, k) - F(n,\, k) =
G(n,\, k + 1) - G(n,\, k),
\end{eqnarray}
where $G(n,\, k) = R(n,\,k)\,F(n,\, k)$ and $\displaystyle
R(n,\,k) =  \frac {2\,k\,(2\,k - 1)}{(n + 2 - k)\,( 2\,n + 5)}$ \
.

By summing both sides of equation (\ref{Thm2eqn1}) with respect to
$k$ we get $S(n + 1)  - S(n) = 0$. Moreover,  $S(-1) = 1$ and
hence $S(n) = 1$ for all $n \in \mathbb{N}_{-1}$ \ .

\bigskip

\textbf{Proof of Theorem \ref{Thm3}.}

Reversing the order of summation, the identity can be rewritten as
\begin{eqnarray}
\label{Thm3eq1} \sum_{m = 2}^{2\,n}\sum_{k = 0}^{\lfloor \frac{m -
2}{2} \rfloor} {{n + k + 1} \choose {m}}2^{m - 1}\,(-1)^{m + k + n
+ 1} =  n\,(n + 1) \ .
\end{eqnarray}

 Let us denote the left side of (\ref{Thm3eq1}) by  $S(n)$ and its summand by $F(n,\, k,\,
m)$, i.e.
\[ F(n,\, k,\, m) = {{n + k + 1} \choose {m}}\,2^{m - 1}\,(-1)^{m + k + n + 1} \ .\]
 Then $F$ satisfies the
recurrence equation:\footnote{The recurrence equation is
automatically generated by  \texttt{MultiSum}, a Mathematica
package  which is available  from
\texttt{htpp//www.risc.uni-linz.ac.at/research/combinat/risc/software/}}
\begin{eqnarray}
\label{Thm3eq2} F(n + 1,\,k,\, m) - F(n,\,k,\,m) =
 F(n,\,k + 1,\, m) - F(n,\,k,\,m) \ .
 \end{eqnarray}

Summing both sides of (\ref{Thm3eq2}) with respect to $k$ and with
respect to $m$, we get
\begin{eqnarray}
\nonumber\label{Thm3eq3} \lefteqn{ S(n + 1)  -
 S(n)} \\
 &  = & \sum_{m = 2}^{2\,n} {{n + \lfloor m/2 \rfloor  + 1}
\choose {m}}\,2^{m - 1}\,(-1)^{m + \lfloor m/2 \rfloor + n + 1}
  - \sum_{m = 2}^{2\,n} {{n + 1} \choose {m}} 2^{m - 1}(-1)^{m + n +
1}.
\end{eqnarray}

But
\begin{eqnarray}
\nonumber\label{Thm3eq4} \sum_{m = 2}^{2\,n} {{n + 1} \choose
{m}}\, 2^{m - 1}\,(-1)^{m + n + 1} & = & \frac{(-1)^{n +
1}}{2}\sum_{m = 0}^{n + 1} {{n + 1} \choose {m}} (-2)^{m} - (n +
1)\,(-1)^{n} + \frac{(-1)^{n}}{2}\\\nonumber
 & = & \frac{(-1)^{n + 1}}{2}(1 - 2)^{n + 1}  - (n +
1)\,(-1)^{n} + \frac{(-1)^{n}}{2}\\
 & = & \frac{1 + (-1)^{n}}{2} - (n + 1)\,(-1)^{n} \ ,
\end{eqnarray}

\begin{eqnarray}
\nonumber\label{Thm3eq5} \lefteqn{\sum_{m = 2}^{2\,n} {{n +
\lfloor m/2 \rfloor  + 1} \choose {m}}\, 2^{m - 1}\,(-1)^{m +
\lfloor m/2 \rfloor + n + 1}}\\\nonumber & = &
 \sum_{m = 1}^{n} {{n +  m  + 1} \choose {2\,m}}\,2^{2\,m - 1}\,(-1)^{m + n + 1}  + \sum_{m =
1}^{n - 1} {{n + m + 1} \choose {2\,m + 1}} \,2^{2\,m}\,(-1)^{m +
n}\\\nonumber & = & \sum_{m = 0}^{n + 1} {{n + m + 1} \choose
{2\,m}}\,2^{2\,m - 1}\,(-1)^{m + n + 1} - 2^{2\,n + 1} +
\frac{(-1)^{n}}{2} \\\nonumber &&+ \sum_{m = 0}^{n} {{n + m + 1}
\choose {2\,m + 1}} \,2^{2\,m}\,(-1)^{m + n}  - (n + 1)\,(-1)^{n}
- 2^{2\,n}\\
 & = & n + \frac{3}{2} - 2^{2\,n + 1} + \frac{(-1)^{n}}{2} + (n +
 1)
 - (n + 1)\,(-1)^{n} - 2^{2\,n} \ .
\end{eqnarray}

From equations (\ref{Thm3eq3}), (\ref{Thm3eq4}) and
(\ref{Thm3eq5}), we get $S(n + 1) -  S(n)  =  2\,(n + 1)$. Since
$S(1) = 2$, and $n\,( n + 1)$ satisfies the same recurrence
relation, therefore  $S(n) = n\,(n + 1)$ for all $n \in
\mathbb{N}_{1}$.

\bigskip

\section*{Some Corollaries.}

For completeness, we state the following results from
(\cite{GZ03}) and prove using theorems 1-3.

\begin{cor}
\label{Cor1}
\[
\sum_{k = 0}^{n}\sum_{m = 2\,k + 1}^{n + k + 1} (-1)^{m + k + n +
1} { {n + k + 1} \choose {m}} 2^{m - 1} = (n + 1)^{2} \ .
\]
\end{cor}

\textbf{Proof:} The result follows by adding theorems 1 and 3.

\begin{cor}
\label{Cor2}

\[
\sum_{k = 0}^{2\,n}\sum_{m = 2\,k + 2}^{2\,n + 2\,k + 2}(-1)^{m +
k} { {2\,n + k + 2} \choose {m}} 2^{m - 1} = (2\,n + 1)(2\,n + 2)
\ .
\]

\end{cor}

\textbf{Proof:} Add theorem \ref{Thm3} and \ref{Cor1} and multiply
the result by 2.

\begin{cor}

\label{Cor3}

\[
\sum_{k = 0}^{2\,n + 1}\sum_{m = 2\,k + 1}^{2\,n + k + 2}(-1)^{m +
k} { {2\,n + k + 2} \choose {m}} 2^{m - 1} = (2\,n + 2)^{2} \ .
\]

\end{cor}

\textbf{Proof:} Replace $n$ by $2\,n + 1$ in theorem \ref{Cor1}.

\begin{cor}

\label{Cor4}

\[
\sum_{k = 0}^{2\,n}\sum_{m = 2\,k + 1}^{2\,n + k + 2}(-1)^{m + k +
1} { {2\,n + k + 1} \choose {m}} 2^{m - 1} = (2\,n + 1)^{2} \ .
\]

\end{cor}

\textbf{Proof:} Replace $n$ by $2\,n$ in theorem \ref{Cor1}.

\begin{cor}
\label{Cor5}

\[
\sum_{k = 0}^{2\,l}\sum_{m =  1}^{2\,k + 2} (-1)^{m}{ {2\,k + m +
1} \choose {2\,m - 1}} 4^{m - 1} = (2\,l + 1)(2\,l + 2) \ .
\]
\end{cor}

\textbf{Proof:} Replace $n$ by $2\,n + 1$ in theorem \ref{Thm1}
and sum $k$ from $0$ to $2\,l$.

\bigskip

{\Large \textbf{Acknowledgement:}} The authors would like to thank
Doron Zeilberger and the referee for their helpful suggestions on
the combinatorial proofs of the main results using the involution
approach.

\end{document}